\documentclass[12pt]{article}
\usepackage[top=2cm,bottom=2cm,
left=2cm,right=2cm]{geometry}

\RequirePackage{amsthm,amsmath,amssymb,dsfont,algorithm,algorithmic}

\RequirePackage[colorlinks,
citecolor=blue,
urlcolor=blue]{hyperref}

\usepackage{bbm,times} 

\theoremstyle{plain}
\newtheorem{thm}{Theorem}[section]
\newtheorem{asm}{Assumption}[section]
\newtheorem{lemma}[thm]{Lemma}
\newtheorem{cor}{Corollary}[section]

\newtheorem{dfn}{Definition}[section]

\newtheorem{rmk}{Remark}[section]

\allowdisplaybreaks

\begin{document}
\title{\bf On high-dimensional classification by sparse generalized Bayesian logistic regression}

\author{The Tien Mai}

\date{
	\begin{small}
		Department of Mathematical Sciences, 
		Norwegian University of Science and Technology,
		\\
		Trondheim 7034, Norway.
		\\
		Email: the.t.mai@ntnu.no
	\end{small}
}

\maketitle

\begin{abstract}
This work addresses the problem of high-dimensional classification by exploring the generalized Bayesian logistic regression method under a sparsity-inducing prior distribution. The method involves utilizing a fractional power of the likelihood resulting the fractional posterior. Our study yields concentration results for the fractional posterior with random design, not only on the joint distribution of the predictor and response variable but also for the regression coefficients. Significantly, we derive novel findings concerning misclassification excess risk bounds using sparse generalized Bayesian logistic regression. These results parallel recent findings for penalized methods in the frequentist literature. Furthermore, we extend our results to the scenario of model misspecification, which is of critical importance.
\end{abstract}


\section{Introduction}

The classification problem plays a central role in statistical learning and has been extensively investigated in various contexts, as exemplified by references such as \cite{devroye1997probabilistic, boucheron2005theory}. Particularly in situations characterized by high-dimensional data, where the number of potential explanatory variables (predictors) $ p $ significantly exceeds the sample size $ n $, a profound challenge arises that transcends disciplinary boundaries, affecting fields such as statistics and machine learning \cite{hastie2009elements, buhlmann2011statistics, fan2010highclassifi, giraud2021introduction}, with various application  as in \cite{chung2010sparse,kotte2020similarity,li2021novel}. Research by \cite{bickel2004some, fan2008high} emphasizes that, even in straightforward cases, high-dimensional classification without feature selection may result in poor performance comparable to random guessing. Consequently,  this issue needs to be addressed by reducing the dimensionality of the feature space through the selective inclusion of a sparse subset of ``meaningful" features. Common approaches in this context often adhere to the frequentist paradigm and center around penalized maximum likelihood methods. Notably, the Lasso and its extensions, detailed in works such as \cite{van2008high, glmnet, hastie2009elements, abramovich2016model}, have been extensively examined and shown to be practically effective. In a more recent study \cite{abramovich2018high}, the authors establish nonasymptotic bounds on misclassification excess risk for procedures based on penalized maximum likelihood.

In the realm of Bayesian literature, research on high-dimensional logistic regression within the Bayesian framework has been carried out. For instance, \cite{jiang2007bayesian, liang2013bayesian,jeong2021posterior,wei2020contraction,atchade2017contraction,narisetty2018skinny} studied fully Bayesian methods in generalized linear models, while \cite{ray2020spike} investigates the variational Bayes method for logistic regression. However, as far as our knowledge extends, there is a notable gap in the literature concerning results on prediction error and misclassification excess risk within the Bayesian framework, comparable to those elucidated in the frequentist literature, as exemplified by \cite{abramovich2018high}. Our work represents an effort to address this gap.

In this study, our primary focus lies in the Generalized Bayesian method, where we employ a fractional power of the likelihood. This gives rise to what is commonly known as fractional posteriors or tempered posteriors, as detailed in \cite{bhattacharya2016bayesian, alquier2020concentration}. Moreover, in these references, fractional posteriors have been shown to provide a more robust approach to address the challenge of model misspecification. It is important to highlight that Generalized Bayesian inference has experienced a surge in attention in recent years, as evidenced by various works such as \cite{matsubara2022robust,hammer2023approximate,jewson2022general,yonekura2023adaptation, medina2022robustness,mai2017pseudo,grunwald2017inconsistency,bissiri2013general,yang2020alpha,lyddon2019general,syring2019calibrating,Knoblauch,mai2023reduced,hong2020model}.

We present important novel findings concerning misclassification excess risk for sparse generalized Bayesian logistic regression in  high-dimensional classification. This result is comparable to those given in the frequentist literature \cite{abramovich2018high}. Additionally, we also investigate the concentration properties of the corresponding fractional posterior. Specifically, we establish concentration results for the joint distribution of the predictor and the response variable in the $\alpha$-Rényi divergence with random design. Consequently, as special cases, we derive concentration results relative to the Hellinger metric and the total variation distance. These concentration results extend on previous work in \cite{jeong2021posterior} in which the authors provide results in Hellinger metric. It also worth to mention that our technical arguments are based on \cite{alquier2020concentration} where they propose a general way to derive concentration rates for fractional posteriors using PAC-Bayesian inequalities.

Diverging from the focus of previous works such as \cite{bhattacharya2016bayesian, alquier2020concentration}, which provided general results for the $\alpha$-Rényi divergence, we extend our investigation to establish concentration results for parameter estimation using specific distance measures. However, this endeavor necessitates the utilization of compatibility numbers, drawn from methodologies in the full Bayesian literature \cite{jeong2021posterior, ray2020spike, castillo2015bayesian}, which have also been previously employed in the frequentist literature \cite{buhlmann2011statistics, van2014higher}. We present precise concentration outcomes regarding $\ell_q$-type metrics directly within the parameter space, where $1 \leq q \leq 2$, these results are comparable to those presented in the frequentist literature \cite{abramovich2016model}. To the best of our knowledge, these results are also innovative within the context of fractional posterior distributions.

In this work, a rescaled Student's t-distribution is employed as the prior, which facilitates inducing sparsity with approximate accuracy. This prior has been utilized in various sparse contexts in previous research \cite{dalalyan2012mirror,dalalyan2012sparse}. This particular prior is advantageous for sampling from the posterior using Langevin Monte Carlo, a gradient-based sampling technique. It is noteworthy that numerous authors have highlighted the significance of heavy-tailed priors in addressing sparsity, as evidenced by several instances \cite{seeger2008bayesian,johnstone2005empirical,rivoirard2006nonlinear,abramovich2007optimality,johnstone2004needles,carvalho2010horseshoe,castillo2012needles,castillo2018empirical,ray2020spike}. Additionally, we briefly demonstrate that a spike and slab prior yields similar theoretical outcomes. However, approximating the corresponding posterior computationally may pose greater challenges.

The remainder of the paper is organized as follows: Section \ref{sc_modelmethod} introduces the high-dimensional classification problem and outlines the generalized sparse Bayesian logistic regression method, along with the associated sparsity-inducing prior distribution. Section \ref{sc_concentration} provides concentration results for the fractional posterior in various metrics and includes extensions to address model misspecification. The findings pertaining to misclassification excess risk are presented in Section \ref{sc_Misclassification}. The paper concludes with a discussion and summary in Section \ref{sc_discuconcluse}. The detailed technical proofs are given in Appendix \ref{sc_proofs}.

\section{Model and method}
\label{sc_modelmethod}
\subsection{Notations and definitions}

We examine a collection of $n$ independent and identically distributed (i.i.d) random variables, represented as $(Z_1, \ldots, Z_n) = Z_1^n$, within a measured sample space $(\mathbb{X}, \mathcal{X}, \mathbb{P})$. The statistical model, denoted by $\left\{P_\theta, \theta \in \Theta\right\}$, encapsulates various probability distributions. The primary goal is to estimate the distribution of the random variables $Z_i$. Initially, we assume the model is well-specified, implying the existence of a parameter value $\theta_0 \in \Theta$ such that $\mathbb{P} \equiv P_{\theta_0}^{\otimes n}$ (specific notifications will be provided later if this assumption is relaxed).

Let $Q$ serve as a dominating measure for this distribution family, and define $p_{\theta} = \frac{{\rm d}P_{\theta}}{{\rm d} Q}(\theta)$. Consider the set $\mathcal{M}_{1}^+(E)$, encompassing all probability distributions on a measurable space $(E,\mathcal{E})$. Assume $\Theta$ is equipped with a $\sigma$-algebra $\mathcal{T}$. Let $\pi\in \mathcal{M}_{1}^+(\Theta)$ be a prior distribution for $ \theta $. The likelihood and the negative log-likelihood ratio will be denoted respectively by
\begin{equation*}
\forall (\theta,\theta')\in\Theta^2\text{, }
L_n(\theta) = \prod_{i=1}^n p_{\theta}(X_i) \text{ and }
r_n(\theta,\theta') = \sum_{i=1}^n \log \frac{p_{\theta'}(X_i)}{p_{\theta}(X_i)}.
\end{equation*}
The fractional posterior, that will be our \textit{ideal} estimator, is given by
\begin{equation}
\label{eq_fractional_posterior}
\pi_{n,\alpha}({\rm d}\theta|X_1^n)
:=
\frac{{\rm e}^{-\alpha r_n(\theta,\theta_0)} \pi({\rm d}\theta)}{\int{\rm e}^{-\alpha r_n(\theta,\theta_0)} \pi({\rm d}\theta)} 
\propto 
L_n^{\alpha}(\theta) \pi({\rm d}\theta),
\end{equation}
as in \cite{bhattacharya2016bayesian,alquier2020concentration}. When $ \alpha =1 $, we recover the usual posterior distribution.

In this work,  a binary regression  model is considered. More specially, let $ Z_i=(Y_i,X_i)\in \lbrace -1,1 \rbrace \times \mathbb{R}^d $ be such that
\begin{equation}
\label{eq_logistic_model}
p_{\theta}
=
\mathbb{P}\lbrace Y=y \vert X=x,\theta\rbrace
=
\frac{e^{yx^\top \theta}}{1+e^{yx^\top\theta}}
.
\end{equation}
We will prove results in the case of random design where we suppose that the distribution of $X_1^n$ does not depend on the parameter. Let $ \|\cdot\|_q $ denote the $ \ell_q $-norm and $ \|\cdot\|_\infty $ denote the max-norm of vectors. We let $ \|\cdot\|_0 $ denote the $ \ell_0 $ (quasi)-norm (the number nonzero entries) of vectors. 

In this work, we study a sparse setting and thus we assume that $ s^* < n $ where $ s^*:= \|\theta_0 \|_0 $.

Let $\alpha\in(0,1)$ and $P,R$ be two probability measures. Let $\mu$ be any measure such that $P\ll \mu$ and $R\ll \mu$. The $\alpha$-R\'enyi divergence  between two probability distributions $P$ and $R$ is  defined by
\begin{align*}
D_{\alpha}(P,R) & =
\frac{1}{\alpha-1} \log \int \left(\frac{{\rm d}P}{{\rm d}\mu}\right)^\alpha \left(\frac{{\rm d}R}{{\rm d}\mu}\right)^{1-\alpha} {\rm d}\mu  \text{,}
\end{align*}
and the Kullback-Leibler (KL) divergence is defined by
\begin{align*}
\mathcal{K}(P,R) & = 
\int \log \left(\frac{{\rm d}P}{{\rm d}R} \right){\rm d}P \text{ if } P \ll R
\text{, }
+ \infty \text{ otherwise}.
\end{align*}

\subsection{Prior specification}

Selecting the appropriate prior distribution plays a pivotal role in achieving a favorable posterior concentration rate in high-dimensional models. In this particular section, we describe a heavy-tailed distribution type that induces the desired concentration rate for the regression coefficient $ \theta $ in high-dimensional settings.
Given a positive number $ C_1 $, for all $ \theta \in B_1 (C_1):= \{ \theta \in \mathbb{R}^d : \|\theta\|_1 \leq C_1 \} $, we consider the following prior,
\begin{eqnarray}
\label{eq_priordsitrbution}
\pi (\theta) 
\propto 
\prod_{i=1}^{d} 
	(\tau^2 + \theta_{i}^2)^{-2}
,
\end{eqnarray}
where $ \tau>0 $ is a tuning parameter.  Technically, it is assumed that $ C_1 > 2d\tau $. This prior has been previously employed in different sparse contexts \cite{dalalyan2012mirror, dalalyan2012sparse}.

It is vital to highlight that $ C_1 $ acts as a regularization constant, generally presumed to be very large. Consequently, the distribution of $ \pi $ closely mirrors that of $ S\tau \sqrt{2} $, where $ S $ is a stochastic vector with independent and identically distributed components originating from the Student's t-distribution with 3 degrees of freedom. By opting for an extremely small $ \tau $, less than $ 1/n $ as in subsequent sections, the bulk of components in $ \tau S $ are drawn near to zero. Nevertheless, due to the heavy-tailed attribute of the Student’s t-distribution, a handful of components of $ \tau S $ are considerably far from zero. This distinct feature equips the prior with the ability to foster sparsity. Furthermore, this prior proves advantageous in facilitating the implementation of Langevin Monte Carlo, a gradient-based sampling technique.

\section{Concentration results}
\label{sc_concentration}
In this section we provide unified results regarding concentration rates of the fractional posterior in high-dimensional logistic regression under suitable assumptions on the design matrix.

\begin{asm}
	\label{asume_finite_2ndmoment}
	  Let assume that $ \mathbb{E}\left\Vert X_1\right\Vert^2 < \infty $.
\end{asm}

\subsection{Results on the distribution}
\label{sc_rslton_distribu}

We first provide our main result concerning the concentration of the fractional posterior relative to the  $\alpha$-R\'enyi divergence of the densities.
  \begin{thm}
\label{thm_main}
Assume that Assumption \ref{asume_finite_2ndmoment} holds. Put $ K_1 := 2\mathbb{E}\left\Vert X_1\right\Vert $, $ \tau = \frac1{n\sqrt{d}} $, and
$$
\varepsilon_n
=
\frac{ K_1 \vee 	4 s^* \log \left(\frac{C_1 n\sqrt{d} }{ s^*}\right) 
	}{n}
	.
$$
For any
    $\alpha\in(0,1)$, $ (\varepsilon,\eta)\in(0,1)^2 $,
 and for all $ \theta_0 $ that $ \|\theta_0\|_1 \leq C_1 - 2d\tau  $,
     we have that
    \begin{equation*}
    \mathbb{P}\left[
 \int D_{\alpha}(P_{\theta},P_{\theta_0}) \pi_{n,\alpha}({\rm d}\theta|X_1^n)  
 \leq  
 \frac{(\alpha+1) \varepsilon_n + \alpha \sqrt{\frac{\varepsilon_n}{n\eta}} + \frac{\log\left(\frac{1}{\varepsilon}\right)}{n}}{1-\alpha}\right]
 \geq 
 1-\varepsilon-\eta.
 \end{equation*}
  \end{thm}

The proof of our main theorem, Theorem \ref{thm_main}, relies on PAC-Bayesian inequalities as argued in \cite{bhattacharya2016bayesian} and \cite{alquier2020concentration}. For a comprehensive overview of this topic, readers may refer to \cite{alquier2024user}. It is worth noting that the concentration of the fractional posterior is achieved solely through the prior concentration rate, see a discussion right after Theorem 2.4 in \cite{alquier2020concentration}. In contrast, the concentration theory of the conventional posterior necessitates more stringent conditions, as highlighted in \cite{jeong2021posterior}. This observation has also been underscored in prior studies such as \cite{bhattacharya2016bayesian}, \cite{alquier2020concentration}, and \cite{cherief2018consistency}.

\begin{rmk}
	It is important to note that the sole assumption regarding the distribution of $X_1$ is that $ K_2:= \mathbb{E}\left\Vert X_1\right\Vert^2 <\infty$. For instance, when $X_1$ follows a uniform distribution on the unit sphere, we observe that $K_1 \leq 2$ and $K_2 \leq 4$. When $X_1 \sim \mathcal{N}(0,s^2 I_d)$, the values are $K_2 = 4 s^2 d$ and $K_1 \leq 2 \sqrt{s^2d}$.
\end{rmk}

Choosing $\eta=\frac{1}{n \varepsilon_n} $ and $\varepsilon = \exp(-n \varepsilon_n)$, we obtain a more readable concentration result by noting that $ 	1-\frac{1}{n\varepsilon_n}-\exp(-n\varepsilon_n) \geq 
1-\frac{2}{n\varepsilon_n} $. Therefore, the sequence $\varepsilon_n$ gives a concentration rate for the fractional posterior in \eqref{eq_fractional_posterior}, stated in the following corollary. 
\begin{cor}
	\label{cor_concentration}
	Under the same assumptions as in Theorem~\ref{thm_main},
	\begin{align*}
	\mathbb{P}\left[
	\int D_{\alpha}(P_{\theta},P_{\theta_0}) \pi_{n,\alpha}({\rm d}\theta|X_1^n)  
	\leq 
	\frac{2(\alpha+1)}{1-\alpha} \varepsilon_n\right] 
	\geq 
	1-\frac{2}{n\varepsilon_n},
	\end{align*}
\end{cor}

\begin{rmk}
	We remind that all technical proofs are given in Appendix \ref{sc_proofs}. Our results imply that the concentration rates are adaptive to the unknown sparsity level $ s^* $.
	While prior studies, as in \cite{jiang2007bayesian, liang2013bayesian, ray2020spike, jeong2021posterior}, concentrate on contraction rates based on the Hellinger distance, our contribution stands out by providing comprehensive concentration results for the fractional posterior using the $\alpha$-Rényi divergence. As specific instances, we deduce the following noteworthy corollary on the Hellinger distance and the total variation distance by leveraging results from \cite{van2014renyi}.
\end{rmk}

Put
\begin{equation}
\label{eq_constantalphaforHellinger}
\mathcal{H}_\alpha 
=
\begin{cases}
\frac{2(\alpha+1)}{1-\alpha}, \alpha \in [0.5,1) ,
\\
\frac{2(\alpha+1)}{\alpha}, \alpha \in (0, 0.5)
.
\end{cases}
\end{equation}
\begin{cor}
	\label{cor_concentration_Hellinger}
	As a special case, Theorem \ref{thm_main} leads to a concentration result in terms of the classical Hellinger distance
	\begin{equation}
	\label{eq_Hellinger_mainresults}
	\mathbb{P}\left[
	\int H^2(P_{\theta},P_{\theta_0}) \pi_{n,\alpha}({\rm d}\theta|X_1^n)  
	\leq 
	\mathcal{H}_\alpha 
	\varepsilon_n \right]
	\geq 
	1-\frac{2}{n\varepsilon_n}
	.
	\end{equation}
	And for $ \alpha \in (0,1) $, 
	\begin{equation}
	\label{eq_dTV_results}
	\mathbb{P}\left[
	\int  d^{2}_{TV}(P_{\theta},P_{\theta_0}) \pi_{n,\alpha}({\rm d}\theta|X_1^n)  
	\leq 
	\frac{4(\alpha+1)}{(1-\alpha)\alpha}
	\varepsilon_n \right]
	\geq 
	1-\frac{2}{n\varepsilon_n}
	,
	\end{equation}
 with $d_{TV}$ being the total variation distance.
\end{cor}

Corollary \ref{cor_concentration_Hellinger} shows that the fractional posterior distribution of $ \theta $ concentrates around its true value at a specified rate relative to the squared Hellinger metric and to the total variation distance. As follow with previous works on factional posterior \cite{alquier2020concentration,bhattacharya2016bayesian}, our results do not require that the true value can be tested against sufficiently separated other values in some suitable sieve, see e.g. \cite{jeong2021posterior}. Nevertheless, the implications of our findings differ to the most current work by \cite{jeong2021posterior}: we deal with random design but for i.i.d observation while \cite{jeong2021posterior} worked with independent observation for fixed design. Moreover, as in \cite{jiang2007bayesian} we assume that $ \|\theta_0\|_1 $ is upper-bounded while \cite{jeong2021posterior} does not assumed.

\begin{rmk}
	Compare to the most recent result from \cite{jeong2021posterior} where they proved result in a deterministic setting, in Hellinger distance, our result is obtained for the case of random design. In the case of high-dimensional setting, $ d>n $, typically our rate is of order $ s^* \log (d /s^* )/n $, while \cite{jeong2021posterior} obtain a rate of order $  s^* \log (d)/n $.
\end{rmk}

We now present a result in expectation, commonly referred to as the consistency result of the fractional posterior.

\begin{thm}
	\label{thm_result_dis_expectation}
For any $\alpha\in(0,1)$, 
	under the same assumptions as in Theorem~\ref{thm_main}, then we have
	\begin{equation*}
	\mathbb{E} \left[ \int D_{\alpha}(P_{\theta},P_{\theta_0}) \pi_{n,\alpha}({\rm d}\theta|X_1^n) \right]
	\leq \frac{1+\alpha}{1-\alpha}\varepsilon_n
	.
	\end{equation*}
\end{thm}

\begin{rmk}
It is noteworthy that, akin to the findings in Corollary \ref{cor_concentration_Hellinger}, one can readily deduce results in expectation for the Hellinger distance or the total variation from Theorem \ref{thm_result_dis_expectation}. For example, we have that $ 	\mathbb{E} \left[ \int H^2(P_{\theta},P_{\theta_0}) \pi_{n,\alpha}({\rm d}\theta|X_1^n) \right]
\leq  \mathcal{H}_\alpha \varepsilon_n $.
\end{rmk}

\subsection{Results on the parameter}
Due to the vagueness associated with the squared Hellinger metric utilized in Corollary \ref{cor_concentration_Hellinger}, no assertion is made regarding the proximity of $ \theta $ and $ \theta_0 $ in the context of a Euclidean-type distance. Diverging significantly from the focus of \cite{jiang2007bayesian,liang2013bayesian}, which exclusively presented results based on the Hellinger metric, our objective is to establish concentration rates for $ \theta $ using a more explicitly defined metric. However, undertaking this task requires the introduction of additional boundedness assumptions.
\begin{asm}
	\label{asume_boundfor_randomdesig}
Assume that there exists a positive constant $ C_0 < \infty $ such that $ |X^\top \theta_0 | < C_0 $, or equivalently there exists  $ 0 < \delta < 1/2 $ such that $ \delta < p_{\theta_0} < 1-\delta $.
\end{asm}

The above assumption has been used before in the context of logistic regression in \cite{abramovich2018high}. Put  $ G:= \mathbb{E}(XX^\top) $.

\begin{thm}
	\label{thm_estimation}
	Under the same assumptions as in Theorem~\ref{thm_main} and additional assume that Assumption \ref{asume_boundfor_randomdesig} holds, for 
	$ 0 \leq \alpha < 1 $, we have that
	\begin{equation}
	\label{eq_estimation_Xbbounds}
\mathbb{E} \left[ \int 
| (\theta - \theta_0 )^\top G  (\theta - \theta_0 ) |
\pi_{n,\alpha}({\rm d}\theta|X_1^n) \right]
\leq  
c \mathcal{H}_\alpha \varepsilon_n 
	,
	\end{equation}
	where $ c $ is a positive universal constant. 
\end{thm}

\begin{rmk}
	The result from Theorem \ref{thm_estimation} is novel. It provides a consistency result in term of an $ \ell_2 $ norm. It is worth mention that this result is different to the most current one from \cite{jeong2021posterior}. Our result is obtained for random design, while \cite{jeong2021posterior} considered fixed designs and provided result in a scaled term as $ \| W_0X^\top (\theta - \theta_0 ) \|^2_2 $, where $ W_0 $ is a diagonal matrix with diagonal elements depending on $ X^\top \theta_0 $.
\end{rmk}

In order to obtain further results in  $ \ell_2 $ distance between, we need to impose additional restriction on the design matrix.

\begin{asm}
	\label{asume_minimaleigenvalue_estimebound}
	Assume that all $ X_{\cdot j} $ are linearly independent. Therefore, the minimal eigenvalue of the matrix $ G:= \mathbb{E}(XX^\top) $, denoted by $ \lambda_{min}(G) $, is strictly positive.
\end{asm}

The above assumption has also been used before in \cite{abramovich2018high}. Under Assumption \ref{asume_minimaleigenvalue_estimebound}, it is straightforward to obtain the following results from Theorem \ref{thm_estimation}.

\begin{cor}
	\label{thm_estimation_ell_1_ell_2}
	Under the same assumptions as in Theorem~\ref{thm_main} and additional assume that Assumption \ref{asume_boundfor_randomdesig} and \ref{asume_minimaleigenvalue_estimebound} holds, for $ 0 \leq \alpha < 1 $, we have
	\begin{equation}
	\label{eq_estimation_ell2bounds}
	\mathbb{E} \left[ \int 
	\| \theta - \theta_0 \|^2_2
	\pi_{n,\alpha}({\rm d}\theta|X_1^n) \right]
	\leq  
	 \frac{c}{\lambda_{min}(G)}
\mathcal{H}_\alpha 
	\varepsilon_n 
	,
	\end{equation}
		where $ c $ is a positive universal constant. 
\end{cor}

\begin{rmk}
In Corollary \ref{thm_estimation_ell_1_ell_2}, our finding regarding the $ \ell_2 $-norm under random design is novel to the best of our knowledge. When compared to the results presented in Theorem 3 of \cite{jeong2021posterior}, which are presented to fixed design scenarios, our derived rate of $ s^*\log (d/s^*)/n $ exhibits a slight improvement over theirs, which stands at $ s^*\log (d)/n $ \cite[Theorem 3]{jeong2021posterior}. It is important to note that their outcomes crucially depend on the concept of the compatibility numbers.
\end{rmk}

\begin{rmk}
It is noteworthy to highlight that within a fixed design framework, the assumptions outlined in Assumption \ref{asume_boundfor_randomdesig} and \ref{asume_minimaleigenvalue_estimebound} can be substituted with broader conditions. An examination of the proof presented in Theorem 3 of \cite{jeong2021posterior} demonstrates the feasibility of incorporating their methodologies into our findings from Theorem \ref{thm_result_dis_expectation}, particularly concerning the evaluation using Hellinger distance. Specifically, we can leverage Lemma A1 from \cite{jeong2021posterior} and utilize the compatibility numbers
$$ \phi_1 (s,W_0) 
:= 
\inf_{\theta:0 < \|\theta\|_0 \leq s} 
\frac{ \| W_0 X^\top\theta \|_2^2 \|\theta\|_0 }{ \|\theta \|_1^2 }
,
\quad
\phi_2 (s,W_0) 
:= 
\inf_{\theta:0 < \|\theta\|_0 \leq s} 
\frac{ \| W_0 X^\top\theta \|_2^2 }{  \|\theta \|_2^2 }
$$
to derive outcomes pertaining to $\ell_1$ and $\ell_2$ distances,  where $ W_0 $ is a diagonal matrix with diagonal elements depending on $ X^\top \theta_0 $.
	
The concept of compatibility numbers finds common usage within Bayesian high-dimensional literature, as highlighted in works such as \cite{castillo2015bayesian, martin2017empirical, ray2020spike, belitser2020empirical}. Originating from high-dimensional frequentist literature, they were initially employed in contexts such as those discussed in \cite{buhlmann2011statistics}.
\end{rmk}

\subsection*{Result in the misspecified case}

In this section, we show that our results can be extended to the misspecified setting.
Assume that the true data generating distribution is parametrized by $\theta_0\notin\Theta$ and define $P_{\theta_0}$ as the true distribution. 
Put
\begin{align*}
\theta^* 
:=
\arg\min_{\theta \in \Theta} 
\mathcal{K}(P_{\theta_0},P_{\theta}),
\end{align*}
we obtain the following result.

\begin{thm}
	\label{thm_misspecified}
For any $\alpha\in(0,1)$, let assume that Assumption \ref{asume_finite_2ndmoment} holds, $ \| \theta^* \|_1 \leq C_1 - 2d\tau $, and with $ \tau = (n\sqrt{d})^{-1} $. Then, 
	\begin{equation*}
	\mathbb{E} \left[ \int D_{\alpha}(P_{\theta},P_{\theta_0}) \pi_{n,\alpha}({\rm d}\theta|X_1^n)\right]
	\leq 
	\frac{\alpha}{1-\alpha} \min_{\theta\in\Theta} \mathcal{K}(P_{\theta_0},P_{\theta})
	+ \frac{1+\alpha}{1-\alpha} r_n
	,
	\end{equation*}
	where 
	$$ 
	r_n = 	\frac{K_1}n
	\vee 
	\frac{	4 \| \theta^* \|_0 \log \left(\frac{C_1 n\sqrt{d} }{\| \theta^* \|_0 }\right)	+\log(2)}{n} 
	.
	$$
\end{thm}

In the well-specified case, i.e $\theta_0 = \arg\min_{\theta\in\Theta} \mathcal{K}(P_{\theta_0},P_{\theta})$, we recover Theorem~\ref{thm_result_dis_expectation}. Otherwise, this result takes the form of an oracle inequality. Although it does not constitute a sharp oracle inequality because the risk measures are not identical on both sides, this observation remains valuable, particularly when $\mathcal{K}(P_{\theta_0},P_{\theta^*})$ is minimal.

In the case of fixed design, however, we can actually derive an oracle inequality result with $ \ell_2 $ error in both sides. The result is as follow.

\begin{cor}
	\label{cor_estimation_miss}
	Under the same assumptions as in Theorem~\ref{thm_misspecified} and additionally assume that Assumption \ref{asume_boundfor_randomdesig} and \ref{asume_minimaleigenvalue_estimebound} hold, for 
$ 0 \leq \alpha < 1 $, we have for the case of fixed design that
	\begin{equation*}
\mathbb{E} \left[ 
	\int \| \theta - \theta_0  \|_2^2
\,
\pi_{n,\alpha}({\rm d}\theta|X_1^n)  
	\right]
	\leq 
	\frac{cK_\alpha}{\lambda_{min}(G) } 
	\min_{\theta\in\Theta} 	
 \| \theta - \theta_0  \|^2_2
	+ 
	\frac{1+\alpha}{1-\alpha} r_n
	,
	\end{equation*}
	where $ c >0 $ is a universal constant.
\end{cor}

Furthermore, put 
\begin{align*}
\hat\theta^B
:=
\int_\Theta \theta \pi_{n,\alpha}({\rm d}\theta|X_1^n)  
\end{align*}
as the mean estimator. By using an
application of Jensen's inequality,
Corollary \ref{cor_estimation_miss} immediately implies an oracle type inequality for the mean estimator,
	\begin{equation*}
\mathbb{E}  
[ \| \hat\theta^B - \theta_0  \|_2^2 ]
\leq 
\frac{cK_\alpha}{\lambda_{min}(G) } 
\min_{\theta\in\Theta} 	
\| \theta - \theta_0  \|^2_2
+ 
\frac{1+\alpha}{1-\alpha} r_n
.
\end{equation*}

\section{Misclassification excess risk bounds}
\label{sc_Misclassification}
In this section, we present novel results regarding the prediction error by using generalized Bayesian logistic regression method. 

We formally consider the following classifier $ \eta_\theta $, for model \eqref{eq_logistic_model},
$$
Y| x 
=
\begin{cases}
1, &\text{ with probability  } p_{\theta}( x),
\\
-1, &\text{ with probability  } 1 - p_{\theta}( x)
\end{cases}
.
$$
The accuracy of a classifier $\eta$
is defined by a misclassification error 
$$
R(\eta)
=
\mathbb{P} (Y \neq \eta( x)).
$$
It is well-known that $R(\eta)$
is minimized by the Bayes classifier 
$ \eta^*( x)= {\rm sign} (p( x) - 1/2) $ \cite{vapnik,devroye1997probabilistic}, i.e.
$
R(\eta^*)
=
\inf R(\eta)
.
$

However, the probability function $p( x)$ is unknown and the resulting classifier $\tilde{\eta}( x)$ should be designed from the data $D_n $: a random sample of $n$ independent observations $( x_1,Y_1),\ldots, ( x_n,Y_n) $, with $ n<d $. 
The design points $ x_i$ may be considered as fixed or random.  
The corresponding (conditional) misclassification error of $\tilde{\eta}$ is given as
$
R(\tilde{\eta})
=
\mathbb{P} (Y \neq \tilde{\eta}( x) \, |D_n)
$
and the goodness of $\tilde{\eta}$ w.r.t. $\eta^*$ is measured by the misclassification
excess risk,
$$
\mathcal{E}(\tilde{\eta},\eta^*)
=
\mathbb{E} \, R(\tilde{\eta})-R(\eta^*).
$$
For logistic regression as in model \eqref{eq_logistic_model}, where it is assumed that
$ p( x) = 1 /( 1+ e^{-\beta^\top x}) $ and $\beta \in \mathbb{R}^d$ is a vector of unknown regression coefficients. The corresponding Bayes classifier 
is a linear classifier that
$
\eta^*( x)
=
{\rm sign} (\beta^{*\top}  x )
.
$
One then estimates $\beta^* $ from the data to get $\hat{\beta}$ (e.g. using maximum likelihood), and
the resulting linear classifier is $\hat{\eta}_{\hat{\beta}} ( x)=  {\rm sign} ( \hat{\beta}^{\top}  x ) $, see e.g. \cite{abramovich2018high}. 

The primary challenges encountered by any classifier are concentrated in the vicinity of the boundary $\{x: p(x)=1/2\}$, which equivalently corresponds to a hyperplane $\beta^\top x=0$ in the logistic regression model. In this region, accurate prediction of the class label becomes particularly challenging. However, in areas where $p(x)$ maintains a substantial distance from $1/2$ (referred to as the margin or low-noise condition), there exists the potential for enhanced bounds on misclassification excess risk.
The following low-noise condition is often made in the classification literature, see e.g. \cite{abramovich2018high,tsybakov2004optimal,mammen1999smooth}.
\begin{asm}
	\label{Assume_low_noise_random}
	Assume that there exist $C>0$ and $\gamma \geq 0$  such that
	\begin{equation}
	\label{eq_low_noise_random}
	P\left( | p_{\theta} -1/2 |
	\leq h\right) 
	\leq 
	C h^{\gamma}
	\end{equation}
	for all $0 < h < h^*$, where $h^* < 1/2$.  
\end{asm}

We are now ready to state our result regarding the misclassification excess risk for the sparse generalized Bayesian method.

\begin{thm}
	\label{thm_missclasification}
	For any $\alpha\in(0,1)$, under Assumption \ref{Assume_low_noise_random} and under the same assumptions as in Theorem \ref{thm_result_dis_expectation}, we have that
	\begin{equation*}
	\int 	
	\mathcal{E} ( \eta_\theta,\eta^*)
	\pi_{n,\alpha}({\rm d}\theta|X_1^n) 
	\leq
C_\alpha 
	\varepsilon_n^{\frac{\gamma + 1}{\gamma+2}}
	,
	\end{equation*}
	where $ C_\alpha>0 $ is a numerical constant depending only on $ \alpha $. 
\end{thm}

\begin{rmk}
	The rate $ 	\varepsilon_n^{(\gamma + 1)/(\gamma+2)} $ is similar to the result in Theorem 7 in \cite{abramovich2018high}. To the best of our knowledge, this result is new for generalized Bayesian logistic regression in high-dimensional sparse classification context. 
\end{rmk}

\begin{rmk}
When $ \gamma = 0 $ in Assumption \ref{Assume_low_noise_random}, it signifies that there is no constraint imposed on the noise. Setting $ C = 1 $ as an example demonstrates that every probability measure complies with this assumption.
\end{rmk}
We immediately obtain the following result without assuming Assumption \ref{Assume_low_noise_random}.
\begin{cor}
	For any $\alpha\in(0,1)$, 	
	\begin{equation*}
	\int 	
	\mathcal{E} ( \eta_\theta,\eta^*)
	\pi_{n,\alpha}({\rm d}\theta|X_1^n) 
	\leq
	C_\alpha \sqrt{\varepsilon_n}	
	,
	\end{equation*}
	where $ C_\alpha>0 $ is a numerical constant depending only on $ \alpha $. 
\end{cor}

\begin{rmk}
It is important to highlight that our derived bounds for misclassification excess risk do not demand any supplementary assumptions, as evidenced in the outcomes for the distribution presented in Section \ref{sc_rslton_distribu}. In contrast, the acquisition of misclassification excess risk bounds for the logistic Lasso \cite{van2008high,abramovich2018high} and the logistic Slope, as discussed in \cite{abramovich2018high}, necessitates specific additional assumptions concerning the characteristics of $X$, such as the weighted restricted eigenvalue or the restricted minimal eigenvalue. This observation constitutes another noteworthy robust finding for the fractional posterior.
\end{rmk}

\section{Discussion and Conclusion}
\label{sc_discuconcluse}

In this work, we have investigated concentration rates within a sparse high-dimensional logistic regression model with random design settings. Employing a sparsity-inducing prior, our focus was on the fractional posterior, achieved by replacing the likelihood with a fractional power of itself. Additionally, we also obtained some results in model misspecification. Novel findings were presented regarding the misclassification excess risk associated with the fractional posterior, which, to the best of our knowledge, are unprecedented. Our outcomes demonstrate comparability with findings documented in the frequentist literature.

Although the primary focus of this study lies in exploring theoretical properties, we offer a concise discussion on the computational aspects associated with the fractional posterior. 
By leveraging the prior specified in equation \eqref{eq_priordsitrbution}, within the realm of logistic regression, we benefit from employing Langevin Monte Carlo (LMC), a gradient-based sampling technique. LMC methods have been effectively showcased in various sparsity contexts, as evidenced by previous studies such as \cite{dalalyan2012sparse,dalalyan2012mirror}. Moreover, LMC has emerged as a promising sampling approach in high-dimensional Bayesian methodologies, as highlighted in works like \cite{durmus2017nonasymptotic,durmus2019high,dalalyan2017theoretical,dalalyan2019user}.

The utilization of the scaled Student's t-distribution, as examined in the preceding section, demonstrates its favorability; however, it lacks the capability for variable selection. 
For sparse scenarios, one may prefer the spike and slab prior of \cite{mitchell1988bayesian,george1993variable}, for variable selection purpose,
\begin{equation}
\label{eq_spike_slab_prior}
\pi_\xi(\theta) 
= 
\prod_{i=1}^{d} \left[p \phi (\theta_i;0,v_1)
+ (1-p) \phi(\theta_i;0,v_0) \right] 
\end{equation}
with $\xi=(p,v_0,v_1)\in[0,1]\times (\mathbb{R}^+)^2$, and $v_0 \ll v_1$. Here, $ \phi (\cdot ; a,b) $ is the Gaussian density with mean $ a $ and variance $ b $.

Employing the spike and slab prior as described in equation \eqref{eq_spike_slab_prior} can yield analogous outcomes to those obtained using the scaled Student prior outlined in equation \eqref{eq_priordsitrbution}. Theorem \ref{thm_spike_slab}, presented below, offers results akin to those derived from the main Theorem \ref{thm_main}.

  \begin{thm}
	\label{thm_spike_slab}
	Let assume that $ \mathbb{E}\left\Vert X_1\right\Vert^2 < \infty $ and $\|\theta_0\|_2 = 1 $. Using spike and slab prior in \eqref{eq_spike_slab_prior} with $ p = 1-e^{-1/d} $, $ v_0 \leq 1/(2n^2d\log(d)) $. Put $ K_1 := 2\mathbb{E}\left\Vert X_1\right\Vert $ and
	$$
	\varepsilon_n
	=
	\frac{ K_1 \vee s^* \log \left(nd\right) 
	}{n}
	.
	$$
Then, 
	\begin{align*}
\mathbb{P}\left[
\int D_{\alpha}(P_{\theta},P_{\theta_0}) \pi_{n,\alpha}({\rm d}\theta|X_1^n)  
\leq 
c_{v_1}\frac{2(\alpha+1)}{1-\alpha} \varepsilon_n\right] 
\geq 
1-\frac{2}{n\varepsilon_n}
,
\end{align*}
where $ c_{v_1} >0 $ is a constant depending only on $ v_1 $.
\end{thm}
The proof is given in Appendix \ref{sc_proofs}.
However, it is noted that developing efficient algorithms for this prior might pose challenges. In the scenario where $ v_0 \rightarrow 0 $, we revert to a more conventional prior that assigns a point mass at zero for each component. Nonetheless, this results in a fractional posterior comprising a mixture of $2^d$ components that blend Dirac masses and continuous distributions, making it more intricate to compute.

\subsection*{Acknowledgments}
The author acknowledges support from the Centre for Geophysical Forecasting, Norwegian Research Council grant no. 309960, at NTNU. 

\subsection*{Conflicts of interest/Competing interests}
The authors declare no potential conflict of interests.

\appendix
\section{Proofs}
\label{sc_proofs}

We first remind some important properties of the $\alpha$-R\'enyi divergence proven in \cite{van2014renyi}.
\begin{rmk}
	\label{rmk_properties_renyi}
	The measures $P$ and $R$ are mutually singular iff $D_{\alpha}(P,R)=(\frac{1}{\alpha-1})\log(0) = +\infty$.
	We have 
	$
	\lim_{\alpha\rightarrow 1} D_{\alpha}(P,R) 
	= 
	\mathcal{K}(P,R)
	$. For $\alpha\in(0,1]$,
	$$ 
	(\alpha/2) d^{2}_{TV}(P,R) \leq D_{\alpha}(P,R) 
	,
	$$ 
	$d_{TV}$ being the total variation distance.  The $\alpha$-R\'enyi divergences are all equivalent for $0<\alpha<1$ and for $\alpha \leq \beta $, 
	$$ 
	\frac{\alpha}{\beta}\frac{1-\beta}{1-\alpha} D_{\beta} \leq D_\alpha \leq D_\beta 
	.
	$$ 
And
	$$
	H^2(P,R) 
	=
	2[1-\exp(-(1/2)D_{1/2}(P,R))] 
\leq
	D_{1/2}(P,R)
	$$ 
where $ H^2(\cdot,\cdot) $ is the squared Hellinger distance.
\end{rmk}

\subsection{Proof of Section \ref{sc_concentration}}

\begin{proof}[\bf Proof of Theorem~\ref{thm_main}]
\label{proof_thm_main}	
	We can check the hypotheses on the KL between the likelihood terms as required in Theorem 2.4 in \cite{alquier2020concentration}. We apply Theorem 2.4 in \cite{alquier2020concentration} for  $\rho_n$ given in \eqref{eq_specific_distribution}.
		Direct calculation shows that the log-likelihood satisfies that
	\[
	\left\vert\log p_\theta(X_1)-\log p_{\theta^\prime}(X_1) \right\vert 
	\leq 
	2\left\Vert X_1\right\Vert
	\left\Vert\theta-\theta^\prime\right\Vert_2
	.
	\]
	Thus, we have 
	\[
	\mathcal{K}(P_{\theta_0},P_\theta)= \mathbb{E}\left[\log p_{\theta_0}(X)-\log p_\theta(X)\right]
	\leq 
	2\mathbb{E}\left\Vert X_1\right\Vert \Vert\theta-\theta_0\Vert_2
	\]
	and 
	\[
	\mathbb{E}\left[\log^2\frac{p_{\theta_0}}{p_\theta}(X)\right]
	= 
	\mathbb{E}\left[\left(\log p_{\theta_0}(X)-\log p_\theta(X)\right)^2\right]
	\leq 
	4\mathbb{E}\left\Vert X_1\right\Vert^2
	\Vert\theta-\theta_0\Vert_2^2
	.
	\]
	  Let assume that $ K_1 := 2\mathbb{E}\left\Vert X_1\right\Vert $ and $ K_2:=4\mathbb{E}\left\Vert X_1\right\Vert^2 < \infty $. 
	When integrating with respect to $\rho_n$ we have 
	\begin{align*}
	\int\mathcal{K}(P_{\theta_0},P_{\theta})\rho_n(d\theta)
	\leq 
	K_1 2\tau \sqrt{d} ,
	\\
	\int\mathbb{E}\left[\log^2\frac{p_{\theta_0}}{p_\theta}(X)\right]\rho_n(d\theta)
	\leq 
	K_2 4\tau^2 d,
	\end{align*}
	by using Lemma \ref{lema_boundfor_ell2}.
	
	To apply Theorem 2.4 in \cite{alquier2020concentration}, it remains to compute the KL between the approximation of the fractional posterior and the prior. From Lemma \ref{lema_boundof_KL}, we have that 
	\begin{align*}
	\frac1n \mathcal{K}(\rho_n,\pi)
	\leq
	\frac{	4 s^* \log \left(\frac{C_1 }{\tau s^*}\right)
		+
		\log(2)}{n}
	.
	\end{align*}
To derive the rate $\varepsilon_n$ as outlined in Theorem 2.4 of \cite{alquier2020concentration}, we combine the aforementioned bounds. Putting $ \tau = \frac1{n\sqrt{d}}\, $, we can apply 
	it with
	\[
	\varepsilon_n
	=
	\frac{K_1}n
	\vee \frac{K_2}{n^2} 
	\vee 
	\frac{	4 s^* \log \left(\frac{C_1 n\sqrt{d} }{ s^*}\right)	+\log(2)}{n}
	.
	\]
	As we have that $ K_1^2 \leq K_2 $, we deduce that 
	\[
	\varepsilon_n
	=
	\frac{K_1}n
	\vee 
	\frac{	4 s^* \log \left(\frac{C_1 n\sqrt{d} }{ s^*}\right)	+\log(2)}{n}
	.
	\]
	The proof is completed.
	
\end{proof}

\begin{proof}[\bf Proof of Corollary \ref{cor_concentration_Hellinger}]
	From Remark \ref{rmk_properties_renyi}, we have that 
	 $$ 
	 H^2(P_{\theta},P_{\theta_0})
	\leq 
	D_{1/2}(P_{\theta},P_{\theta_0})  
	\leq 
	D_{\alpha}(P_{\theta},P_{\theta_0}) 
	,
	$$
	for $ \alpha \in [0.5,1) $. In addition, we also have that 
	$$ 
	D_{1/2}(P_{\theta},P_{\theta_0})  
	\leq 
	\frac{(1-\alpha)1/2}{\alpha (1-1/2)} D_{\alpha}(P_{\theta},P_{\theta_0}) 
	=
	\frac{(1-\alpha)}{\alpha} D_{\alpha}(P_{\theta},P_{\theta_0}) 
	,
	$$ 
	for $ \alpha \in (0, 0.5)  $.
	
	Thus, using definition of $ \mathcal{H}_\alpha  $ and Corollary \ref{cor_concentration}, we obtain the results.
	
\end{proof}

\begin{proof}[\bf Proof of Theorem~\ref{thm_result_dis_expectation}]
	With similar steps as in the proof of Theorem \ref{thm_main}, but now we apply Theorem 2.6 in \cite{alquier2020concentration}. Therefore, we only need to use the following conditions
	\begin{align*}
	\int\mathcal{K}(P_{\theta_0},P_{\theta})\rho_n(d\theta)
	\leq 
	\varepsilon_n
	,
	\end{align*}
	and
	\begin{align*}
	\mathcal{K}(\rho_n,\pi)
	\leq
	n\varepsilon_n
	.
	\end{align*}
	which can be found in the proof of Theorem \ref{thm_main}. This completes the proof.
	
\end{proof}

\begin{proof}[\bf Proof of Theorem \ref{thm_estimation}]
	Under Assumption \ref{asume_boundfor_randomdesig}, it can be verified that, for any $ \theta_1, \theta_2 $ and the corresponding $ p_{\theta_1} $ and $ p_{\theta_2} $,
	\begin{align*}
	\delta (1-\delta) | X^\top (\theta_1 - \theta_2 ) | 
	\leq 
	| p_{\theta_1} - p_{\theta_2} |,
	\end{align*}
	(details are given in inequality (51) in the proof of Theorem 7 in \cite{abramovich2018high}). The above inequality leads to
	\begin{align*}
	\| p_{\theta_1} - p_{\theta_2} \|_{L_2(q)}
	\geq
	\delta (1-\delta) 
	\sqrt{(\theta_1 - \theta_2 )^\top G  (\theta_1 - \theta_2 )}
	,
	\end{align*}
	where $ \|h\|_{L_2(q)} = (\int h^2(x)q(x)dx)^{1/2} $ is the $ L_2$-norm of $ h $ weighted by the marginal distribution $ q $ of $ X $ (see also inequality (52) in \cite{abramovich2018high}).	Now, from inequality (49) in \cite{abramovich2018high}, we have that
	\begin{align*}
	H^2(P_{\theta_1},P_{\theta_2}) 
	\geq
	\frac{1}{2} 
	\| p_{\theta_1} - p_{\theta_2} \|^2_{L_2(q)}
	.
	\end{align*}
	As a consequence of Theorem \ref{thm_result_dis_expectation}, we have that
	\begin{align*}
	\mathbb{E} \left[ \int H^2(P_{\theta},P_{\theta_0}) \pi_{n,\alpha}({\rm d}\theta|X_1^n) \right]
	\leq  
	\mathcal{H}_\alpha \varepsilon_n 
	\end{align*}
	and thus we obtain that
	\begin{align*}
	\mathbb{E} \left[ \int 
	| (\theta - \theta_0 )^\top G  (\theta - \theta_0 ) |
	\pi_{n,\alpha}({\rm d}\theta|X_1^n) \right]
	\leq  
	\frac{2}{	\delta^2 (1-\delta)^2 }
	\mathcal{H}_\alpha \varepsilon_n 
	.
	\end{align*}
 This completes the proof.
	
\end{proof}

\begin{proof}[\bf Proof of Theorem~\ref{thm_misspecified}]
	
	Direct calculation shows that the log-likelihood satisfies that
	\[
	\left\vert\log p_\theta(X_1)-\log p_{\theta^\prime}(X_1) \right\vert 
	\leq 
	2\left\Vert X_1\right\Vert
	\left\Vert\theta-\theta^\prime\right\Vert_2
	.
	\]
	We can check the assumptions as in Theorem 2.7 in \cite{alquier2020concentration}.
	First, we have that
	\[
	\mathbb{E}_{\theta_0}\left[\log\frac{ p_{\theta^*}(X_i)}{ p_{\theta}(X_i)}\right]
	= 
	\mathbb{E}\left[\log p_{\theta^*}(X)-\log p_\theta(X)\right]
	\leq 
	2\mathbb{E}\left\Vert X_1\right\Vert \Vert\theta-\theta^* \Vert_2
	.
	\]
	Let assume that $ K_1 := 2\mathbb{E}\left\Vert X_1\right\Vert $ and $ K_2:=4\mathbb{E}\left\Vert X_1\right\Vert^2 < \infty $. 
	
	We now consider the following distribution as a translation of the prior $ \pi $,
	\begin{equation*}
	\rho^*_{n} (\theta) 
	\propto 
	\pi (\theta - \theta^*)\mathbbm{1}_{B_1(2d\tau)} (\theta - \theta^*).
	\end{equation*}
	Given $ \| \theta^* \|_1 \leq C_1 - 2d\tau $, as $ \theta - \theta^*\in B_1(2d\tau) $ holds, it implies that $ \theta \in B_1(C_1) $. Consequently, the distribution $ \rho^*_{n} $ is absolutely continuous with respect to the prior distribution $ \pi $. 	
	From Lemma \ref{lema_boundfor_ell2}, 
	\begin{align*}
	\int \mathbb{E}_{\theta_0}\left[\log\frac{ p_{\theta^*}(X_i)}{ p_{\theta}(X_i)}\right]
	\rho^*_n({\rm d}\theta) 
	\leq 
	K_1 2\tau \sqrt{d}
	.
	\end{align*}
	
	To apply Theorem 2.7 in \cite{alquier2020concentration} it remains to compute the KL between the approximation of the fractional posterior and the prior. From Lemma \ref{lema_boundof_KL}, we have that 
	\begin{align*}
	\frac1n \mathcal{K}(\rho_n,\pi)
	\leq
	\frac{	4 \| \theta^* \|_0 \log \left(\frac{C_1 }{\tau  \| \theta^* \|_0}\right)
		+
		\log(2)}{n}
	.
	\end{align*}
	To obtain an estimate of the rate $\varepsilon_n$ of Theorem 2.7 in \cite{alquier2020concentration} we put together those bounds. Choosing $ \tau = \frac1{n\sqrt{d}} \, $, we can apply 
	it with
	\[
	\varepsilon_n
	=
	\frac{K_1}n
	\vee 
	\frac{	4 \| \theta^* \|_0 \log \left(\frac{C_1 n\sqrt{d} }{\| \theta^* \|_0 }\right)	+\log(2)}{n}
	.
	\]
	The proof is completed.
	
\end{proof}

\begin{proof}[\bf Proof of Corollary~\ref{cor_estimation_miss}]
	Using Lemma 1 in \cite{abramovich2016model}, there exist a positive constant $ c>0 $ that
	\[
	\mathcal{K}(P_{\theta_0},P_\theta)
	\leq 
c \| \theta-\theta_0 \|^2_2
	.
	\]
From inequalities (52) and (49) from the proof of Theorem 7 in \cite{abramovich2018high}, we have that
	\[
\frac{	
	\delta (1-\delta) \lambda_{min} (G) }{2} 
\|\theta - \theta_0  \|_2^2 
	\leq
	H^2(P_{\theta},P_{\theta_0})  
	\leq
	\mathcal{H}_\alpha D_{\alpha}(P_{\theta},P_{\theta_0})
	.
	\]
	Thus the result of the corollary is obtained directly from Theorem \ref{thm_misspecified}.
	
\end{proof}

\subsection{Lemmas}

\begin{dfn}
	We define the following distribution as a translation of the prior $ \pi $,
	\begin{equation}
	\label{eq_specific_distribution}
	p_0(\theta) 
	\propto 
	\pi (\theta - \theta_0 )
	\mathbbm{1}_{B_1(2d\tau)} (\theta - \theta_0 ).
	\end{equation}
\end{dfn}
It is worth highlighting that given $ \| \theta_0 \|_1 \leq C_1 - 2d\tau $, when the condition $ \theta - \theta_0 \in B_1(2d\tau) $ holds, it implies that $ \theta \in B_1(C_1) $. Consequently, the distribution $ p_0 $ is absolutely continuous with respect to the prior distribution $ \pi $ and thus the $ KL $ is finite.

The proofs the following two lemmas can be found \cite{mai2023high}. They are derived by using results from \cite{dalalyan2012mirror} (using Lemma 2 and 3).

\begin{lemma}
	\label{lema_boundof_KL}
	Let $p_0$ be the probability measure defined by (\ref{eq_specific_distribution}). Then
	$$
	KL(p_0,\pi)
	\leq
	4 s^* \log \left(\frac{C_1 }{\tau s^*}\right)
	+
	\log(2)
	.
	$$
\end{lemma}
\begin{lemma}
	\label{lema_boundfor_ell2}
	Let $p_0 $ be the probability measure defined by (\ref{eq_specific_distribution}). If
	$d\geq 2$ then
	$$
	\int_\Lambda \| \beta- \beta^* \|^2 p_0(d \beta)
	\leq
	4d\tau^2 
	.
	$$
\end{lemma}

\subsection{Proofs of Section \ref{sc_Misclassification}}

\begin{proof}[\bf Proof of Theorem~\ref{thm_missclasification}]
	For any $\alpha\in(0,1)$, from Theorem \ref{thm_result_dis_expectation},
	\begin{equation*}
	\mathbb{E} \left[ \int D_{\alpha}(P_{\theta},P_{\theta_0}) \pi_{n,\alpha}({\rm d}\theta|X_1^n) \right]
	\leq \frac{1+\alpha}{1-\alpha}\varepsilon_n.
	\end{equation*}
	Then, Fubini's theorem leads to
	\begin{equation*}
	\int \mathbb{E}  D_{\alpha}(P_{\theta},P_{\theta_0}) \pi_{n,\alpha}({\rm d}\theta|X_1^n) 
	\leq 
	\frac{1+\alpha}{1-\alpha}\varepsilon_n.
	\end{equation*}
	Then, 
	\begin{equation*}
	\int \mathbb{E}  
	H^2(P_{\theta},P_{\theta_0})
	\pi_{n,\alpha}({\rm d}\theta|X_1^n) 
	\leq 
\mathcal{H}_\alpha
		\varepsilon_n
		.
	\end{equation*}
	Using inequalities (48) and (49) in \cite{abramovich2018high}, that make use of Theorem 3 in \cite{bartlett2006convexity}, under Assumption \ref{Assume_low_noise_random}, there exists a constant $ C>0 $ such that
	\begin{equation*}
	\mathcal{E}(\hat{\eta},\eta^*)
	\leq
	C \left(\mathbb{E}  
	H^2(P_{\theta},P_{\theta_0})
	\right)^{\frac{\gamma + 1}{\gamma+2}}
	.
	\end{equation*}
	Thus, we obtain the result of the theorem.
	
\end{proof}

\begin{proof}[\bf Proof of Theorem \ref{thm_spike_slab}] Under the assumption that $\|\theta_0\|_2 = 1 $.
	We start by defining, 
	for any  $\delta>0$, 
	$$ 
	\rho_{\theta_0,\delta}
	({\rm d}\theta) 
	\propto
	\mathbf{1}_{\|\theta-\theta_0\|\leq \delta} \pi({\rm d}\theta) 
	.
	$$
	From the proof of Theorem 2.5 in \cite{ridgway2014pac}, we have for $ v_0 \leq \delta^2 /(2d\log(d)) $ that 
	\begin{align*}
	\mathcal{K}(\rho_{\theta_0,\delta},\pi)
	= 
	\|\theta_0\|_0 \left[ \log\left(\frac{2 \sqrt{2 \pi v_1 d}}{ p \delta}\right) + \frac{1}{v_1}
	+ \frac{\delta^2}{v_1 d} \right]+ \log(2) + d\log\frac{1}{1-p}.
	\end{align*}
	Let assume that $ K_1 := 2\mathbb{E}\left\Vert X_1\right\Vert $ and $ K_2:=4\mathbb{E}\left\Vert X_1\right\Vert^2 < \infty $. 
	Similarly to the proof of Theorem \ref{thm_main}, page \pageref{proof_thm_main}, integrating with respect to $ \rho_n := \rho_{\theta_0,\delta} $ we have 
	\begin{align*}
	\int\mathcal{K}(P_{\theta_0},P_{\theta})
\rho_n d\theta)
	& \leq 
	K_1 2 \delta
	,
	\\
	\int\mathbb{E}\left[\log^2\frac{p_{\theta_0}}{p_\theta}(X)\right]
\rho_n(d\theta)
	& \leq 
	K_2 4\delta^2 
	.
	\end{align*}
	Taking $ p = 1-e^{1/d} $, $ v_0 \leq 1/(2n^2d\log(d)) $ and $ \delta = 1/n $, we have that
	\begin{align*}
	\frac1n \mathcal{K}(\rho_n,\pi)
	\leq
	c_{v_1}
	\frac{	 s^* \log \left(nd\right)
	}{n}
	.
	\end{align*}
To calculate  the rate $\varepsilon_n$ as specified in Theorem 2.4 of  \cite{alquier2020concentration}, we combine these bounds. Choosing $ \delta = n^{-1} $, we can apply 
	it with
	\[
	\varepsilon_n
	=
	\frac{K_1}n
	\vee \frac{K_2}{n^2} 
	\vee 
	\frac{	 s^* \log \left(nd\right)	}{n}
	.
	\]
	As we have that $ K_1^2 \leq K_2 $, we deduce that $ \varepsilon_n $ is with
	\[
	\varepsilon_n
	=
	\frac{K_1}n
	\vee 
	\frac{	4 s^* \log \left(nd\right)	}{n}
	.
	\]
	The proof is completed.
\end{proof}

\end{document}